%% file: Dav3.tex
\documentclass{amsart}
\include{Dav.board}

\newcommand{\N}{\mathbb{N}}
\newtheorem{Theo}{Theorem}
\newtheorem{prop}{Proposition}
\newtheorem{Lem}{Lemma}

\begin{document}
\title[Davenport's constant]{Davenport's constant for groups of the form
  $\Z_3\oplus\Z_3\oplus\Z_{3d}$}
\author[G. Bhowmik]{Gautami Bhowmik}
\author[J.-C. Schlage-Puchta]{Jan-Christoph Schlage-Puchta}
\begin{abstract}
We determine Davenport's constant for all groups of the form $\Z_3\oplus
\Z_3\oplus\Z_{3d}$.
\end{abstract}
\maketitle

\section{Introduction and Notation}

For a finite abelian group $G$ let $D(G)$ be the Davenport's constant,
that is, the least integer $n$ such that among each sequence $g_i$ in
$G$ there exists a non-empty subsequence $g_{i_k}$ with sum 0. Writing $G$ as
\[
G \cong \Z_{d_1}\oplus\Z_{d_2}\oplus \dots \oplus \Z_{d_r}, \qquad d_1|d_2|\dots |d_r,
\]
we obtain a sequence of $\sum_i d_i -r$ elements without
a zerosum subsequence, thus we have the trivial bound $D(G)\geq M(G) = \sum_i d_i -r+1$.
It has been conjectured that $D(G)=M(G)$ holds true for all finite
groups, and this conjecture was proven for various special cases,
including finite $p$-groups, and groups of rank $r\leq 2$. However, there
are infinitely many counterexamples known for every rank $r\geq 4$. It is unknown
whether $D(G)=M(G)$ holds true for all groups of rank 3, the authors
are inclined to believe that this is always the case. The simplest
undecided case up to now is $G=\Z_3\oplus\Z_3\oplus\Z_{15}$, which was already
mentioned by van Emde Boas and Kruyswijk \cite{Boas1}.
In the present note this case is solved, more generally, we show the
following. 
\begin{Theo}
\label{thm:D=MZ3}
Let $d$ be an integer, $A\subseteq\Z_3\oplus\Z_3\oplus\Z_{3d}$ be a multiset
consisting of $3d+4$ elements. Then there exists a multiset $B\subseteq A$,
such that $\sum_{b\in B} b=0$. 
\end{Theo}

Our approach is inspired by an idea of Delorme, Ordaz and
Quiroz \cite{DOQ}. Suppose we are given a sequence $A$ of $3d+4$ points in
$\Z_3\oplus\Z_3\oplus\Z_{3d}$. Consider the image $\tilde{A}$ of this sequence under
the canonical projection $\Z_3\oplus\Z_3\oplus\Z_{3d}\to \Z_3^3$. If this
sequence contains a family of $d$ pairwise disjoint subsequences
adding up to zero, we obtain a sequence of $d$ elements in $\Z_d$,
each of which is represented as a sum of certain elements in
$A$. Among these elements we choose a subsequence adding up to 0, and
find that $A$ contains a subsequence adding up to 0. Using this method
Delorme, Ordaz and Quiroz showed that for groups of the form
$G=\Z_3\oplus\Z_3\oplus\Z_{3d}$ we have $D(G)\leq M(G)+2$. Unfortunately,
this inequality is the best possible, since for every $d\geq 3$ there exists
a sequence $A\subseteq G$ with $3d+5$ elements, which does not contain $d$
pairwise disjoint zerosum subsets. To remedy this, we note that for
$(d, 3)=1$ we have 
\[
\Z_3\oplus\Z_3\oplus\Z_{3d} \cong \Z_3^3\oplus \Z_d,
\]
thus we can represent a sequence $A\subseteq\Z_3\oplus\Z_3\oplus\Z_{3d}$ as a pair
$(\tilde{A}, f)$, where $\tilde{A}\subseteq\Z_3^3$ is the image of $A$ under
the canonical projection, and $f:\tilde{A}\to\Z_d$ is the function such
that the element $a_i\in A$ is represented as $(\tilde{a}_i,
f(\tilde{a}_i))$ in $\Z_3^3\oplus \Z_d$. This idea allows us to concentrate
only on the small group $\Z_3^3$; in fact, the remainder of this
article deals only with combinatorial properties of $\Z_3^3$.

Since the order of elements plays no r\^ole, we will henceforth speak
of multisets instead of sequences. To visualize the combinatorial
considerations, we view $\Z_3^3$ as the elements
of a $3\oplus 3\oplus 3$-cube, and this cube again as three $3\oplus 3$-squares
placed side by side. The origin is placed in the lower left corner of
the leftmost 
rectangle and coordinates are associated to cells in the order of board, row and
column. Cells marked with a black circle are
elements that are known to be contained in the set under consideration, cells
marked with a white circle denote elements known not to be contained in the
set. Cells marked with a black circle and a white number $n$ denote elements
which are known to be contained in the multiset under consideration at least
$n$ times. For example, in the following visualization of some information on
a multiset  $A\subseteq\Z_3+\Z_3+\Z_3$, the
cell marked o denotes the neutral element of the group, whereas the cell
marked a is the element $(1, 2, 1)$, and the multiset $A$ contains the element
$(1, 1, 1)$, and the element $(0, 1, 2)$ at least twice.

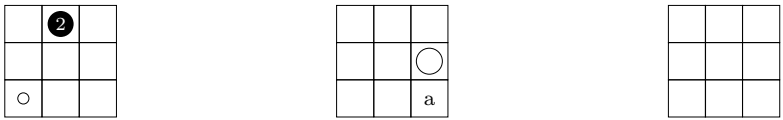
\begin{miniboard}
+ 2 +
+ + +
o + +
\end{miniboard}\hspace*{\fill}
\begin{miniboard}
+ + +
+ + O
+ + a
\end{miniboard}\hspace*{\fill}
\begin{miniboard}
+ + +
+ + +
+ + +
\end{miniboard}

One advantage of this notation is the fact that one can often read off
the existence of zero-sums from the picture. For example, 3 distinct elements
add up to 0 if and only if they lie on an affine line, thus, in the
following picture the three cells marked a as well as the three cells
marked b are zerosum subsets.

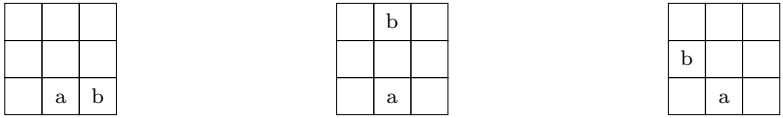
\begin{miniboard}
+ + +
+ + +
+ a b
\end{miniboard}\hspace*{\fill}
\begin{miniboard}
+ b +
+ + +
+ a +
\end{miniboard}\hspace*{\fill}
\begin{miniboard}
+ + +
b + +
+ a +
\end{miniboard}

We will use this argument repeatedly without further explanation.

The technique of  Delorme, Ordaz, and Quiroz requires the study of certain
auxilliary functions. Denote by $D_k(G)$ the least integer $n$, such
that every multiset of $n$ elements contains $k$ disjoint zerosum
subsets, and by $D^k(G)$ the least integer $n$, such that every
multiset of $n$ elements contains a zerosum subset consisting of at
most $k$ elements. Note that $D^k(G)$ is finite only  if $k$ is at least
the exponent of $G$; the case that $k$ equals the exponent of $G$ has
received particular interest.
\footnote{D* already exists in literature. We could use $\bar D$ or Ol}
Adding a star always means that we ask for the least $n$, such that
each subset of $n$ distinct elements has the reuired property, for
example, $ D^*_2(\Z_4)=4$, since 4 distinct points in $\Z_4$ contain
the element 0 as well as the set $\{1, 3\}$, and therefore 2 disjoint
zerosum subsets. In particular, $D^*(G)$ is known as Olson's constant,
which was determined for cyclic groups by Olson \cite{Olson}.

This article is organized as follows. In the next section we give
several rather special results for the variations of $D(\Z_3^3)$ just
mentioned.  These results are of limited interest, but will save us a
lot of work later on. In the last section
we describe the splitting of the set $A$ into the pair $(\tilde{A},
f)$, and prove Theorem~\ref{thm:D=MZ3}.

Our approach is not restricted to groups of the form
$\Z_3\oplus\Z_3\oplus\Z_{3d}$, but can be applied to all sequences of the form
\[
\Z{a_1}\oplus\Z_{a_2}\oplus\dots\oplus\Z_{a_{r-1}}\oplus\Z_{a_rd}, \qquad a_1|a_2|\dots|a_r,
\]
where $a_1, \ldots, a_r$ are fixed, and $d$ runs over all integers coprime
to $a_r$. However, soon the computational effort becomes  too large for
a treatment as explicit as given here. In work in progress, we hope to
automatize parts of the proof to deal with larger groups as well.

\section{Some special values of $D_k$ and related functions}

The results of this section are summarized in the following.
\begin{prop}
We have
\[
\begin{array}{rclrclrcl}
D^3(\Z_3^3) & = & 17 & D^4(\Z_3^3) & = & 10 & D^5(\Z_3^3) & = & 9\\
D^{3*}(\Z_3^3) & = & 9 & D^{4*}(\Z_3^3) & = & 8 & D^{5*}(\Z_3^3) & = & 8\\
D_2(\Z_3^3) & = & 11 & D_k(\Z_3^3) & = & 3k+6 (k\geq 3)\\
D^*(\Z_3^3) & = & 7 & D_2^*(\Z_3^3) & = & 10\\
\end{array}
\]
\end{prop}
\begin{Lem}
\label{Lem:distinctfew}
Set $G=\Z_3+\Z_3+\Z_3$.
\begin{enumerate}
\item  Let $A=\{a_1, \ldots, a_6\}$ be a set of distinct elements of
$G$ such that there does not exist a zerosum subset $Z$ of $A$ with at most
3 elements. Then there are distinct indices $i, j, k$, such that $a_i+a_j=a_k$.
\item Let $A=\{a_1, \ldots, a_8\}$ be a set of distinct elements of
$G$ such that there does not exist a zerosum subset $Z$ of $A$ with at most
3 elements. Then, up to linear equivalence, $A$ is the set

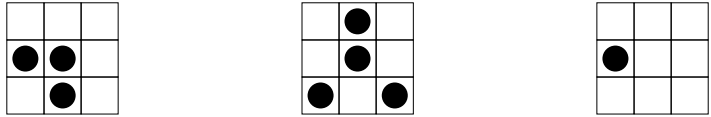
\begin{miniboard}
+ + +
X X +
+ X +
\end{miniboard}\hspace*{\fill}
\begin{miniboard}
+ X +
+ X +
X + X
\end{miniboard}\hspace*{\fill}
\begin{miniboard}
+ + +
X + +
+ + +
\end{miniboard}

\item Let $A=\{a_1, \ldots, a_9\}$ be a set of distint elements of
$G$. Then there exists a zerosum subset $Z$ of $A$ with at most
3 elements.
\end{enumerate}
\end{Lem}
\begin{proof}
Let $ A=\{a_1, \ldots, a_6\}$ be a set of 6 elements, and suppose that none of 
the equations $x+y=z$ and $x+y+z=0$ is solvable within $A$. Then $A$ cannot be
contained in a plane, 
thus, we may choose a basis in $A$. We therefore obtain the following
description of $A$.

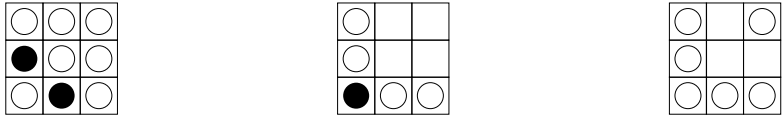
\begin{miniboard}
O O O
X O O
O X O
\end{miniboard}\hspace*{\fill}
\begin{miniboard}
O + +
O + +
X O O
\end{miniboard}\hspace*{\fill}
\begin{miniboard}
O + O
O + +
O O O
\end{miniboard}

Any two of the three points $(1, 1, 2)$, $(1, 2, 1)$ and $(2, 1, 1)$ form
together with one of the points $(0, 0, 1)$, $(0, 1, 0)$ and $(1, 0, 0)$ a
zerosum subset, thus, by symmetry we may assume that $(1, 1, 2)$ and $(1, 2,
1)$ are not contained in $A$. Moreover, if $(2, 1, 1)$ were in $A$, the only
remaining position would be $(1, 2, 2)$, and we would have $|A|\leq 5$, that
is, $(2, 1, 1)$ is not in $A$ as well. From the remaining 4 positions, 3 have
to be taken by elements in $A$, but $(1, 1, 1)$ and $(1, 2, 2)$ cannot be
taken at the same time, thus, both $(2, 1, 2)$ and $(2, 2, 1)$ have to be in
$A$. Then $(1, 1, 1)$ cannot be contained in $A$, and we obtain the following
situation. 

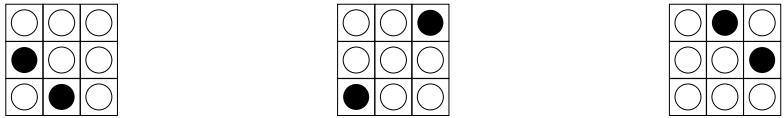
\begin{miniboard}
O O O
X O O
O X O
\end{miniboard}\hspace*{\fill}
\begin{miniboard}
O O X
O O O
X O O
\end{miniboard}\hspace*{\fill}
\begin{miniboard}
O X O
O O X
O O O
\end{miniboard}

But now we have $(2, 1, 2)+(1, 2, 2)=(0, 0, 1)$, proving our claim.

Now let $A$ be a set of size 8 without a zerosum of length 3. By part 1 we may
assume that $a_1+a_2=a_3$, moreover, not all elements of $A$ are contained in
the plane generated by $a_1$ and $a_2$, and we obtain the following situation.

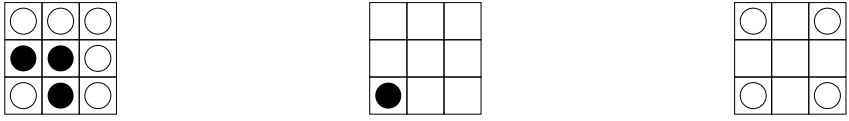
\begin{miniboard}
O O O
X X O
O X O
\end{miniboard}\hspace*{\fill}
\begin{miniboard}
+ + +
+ + +
X + +
\end{miniboard}\hspace*{\fill}
\begin{miniboard}
O + O
+ + +
O + O
\end{miniboard}
Moreover, we may suppose that there are more elements in the middle layer than
in the uppermost one. Suppose first that $(1, 1, 1)$ is in $A$. 

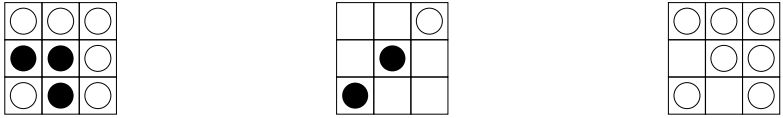
\begin{miniboard}
O O O
X X O
O X O
\end{miniboard}\hspace*{\fill}
\begin{miniboard}
+ + O
+ X +
X + +
\end{miniboard}\hspace*{\fill}
\begin{miniboard}
O O O
+ O O
O + O
\end{miniboard}

If both the remaining cells in the uppermost layer were contained in $A$, then
no further cell in the middle layer could be contained in $A$; if on the other
hand both cells were not contained in $A$, then there are three more cells in
the midle layer, and it is easly seen that this implies the existence of a
zerosum sequence of length 3 in the middle layer. Hence, precisely one of $(2,
0 ,1)$ and $(2, 1, 0)$ is in $A$, and we may assume that this element is $(2,
0, 1)$. Then we reach the following situation.

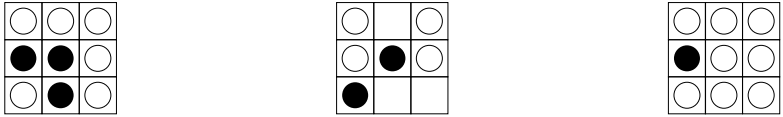
\begin{miniboard}
O O O
X X O
O X O
\end{miniboard}\hspace*{\fill}
\begin{miniboard}
O + O
O X O
X + +
\end{miniboard}\hspace*{\fill}
\begin{miniboard}
O O O
X O O
O O O
\end{miniboard}

From the remaining three cells two have to be taken, but $(1, 1, 0)$ would
yield with one of the other two cells and $(1, 0, 0)$ resp. $(1, 1, 1)$ a
zerosum, thus, we obtain the constellation given in the Lemma.

Now suppose that $(1, 1, 1)$ is not in $A$. Since any element in the two upper
layers can be interchanged with $(1, 0, 0)$ by a linear transformation leaving
the lower layer fixed, we can avoid this case unless for all elements $x, y,
z\in A$ with $x+y\in A$ and $z\neq x+y$ we have $x+y+z\not\in A$. In
particular, in our situation this implies that we may suppose that for each
$z\in A$ which is not in the lower layer both elements $z\pm(0, 1, 1)$ are not
in $A$. Assume that $(1, 1, 0)$ is in $A$. The

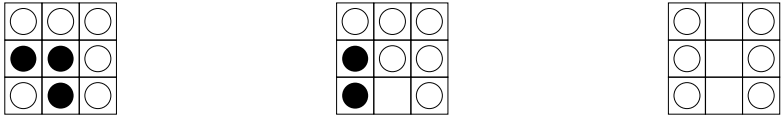
\begin{miniboard}
O O O
X X O
O X O
\end{miniboard}\hspace*{\fill}
\begin{miniboard}
O O O
X O O
X + O
\end{miniboard}\hspace*{\fill}
\begin{miniboard}
O + O
O + O
O + O
\end{miniboard}

Sine there are at least 3 elements in the middle layer, we deduce that $(1, 1,
0)\in A$, which contradicts the fact that there are two more elements in the
uppermost layer. Next suppose that $(1 ,0, 2)$ was in $A$. Then we obtain the
following. 

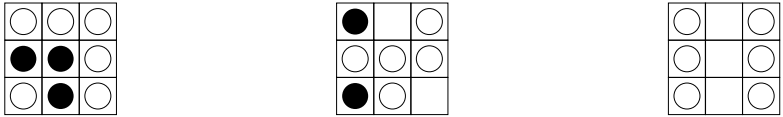
\begin{miniboard}
O O O
X X O
O X O
\end{miniboard}\hspace*{\fill}
\begin{miniboard}
X + O
O O O
X O +
\end{miniboard}\hspace*{\fill}
\begin{miniboard}
O + O
O + O
O + O
\end{miniboard}

In the middle layer there has to be another element of $A$, however ,for both
possible places we see that there could be at most one other cell in the
uppermost layer. Hence, $(1, 0, 2)$ and $(1, 2, 0)$ are not in $A$.
Then we obtain the following constellation. 

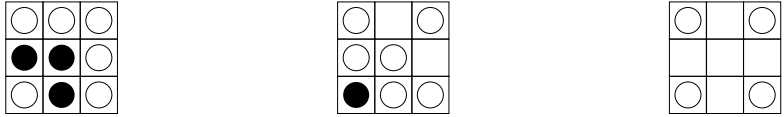
\begin{miniboard}
O O O
X X O
O X O
\end{miniboard}\hspace*{\fill}
\begin{miniboard}
O + O
O O +
X O O
\end{miniboard}\hspace*{\fill}
\begin{miniboard}
O + O
+ + +
O + O
\end{miniboard}

Here, not both cells in the middle
layer can be taken, thus there have to be three cells in the uppermost layer,
contrary to our assumption that there are more cells in the middle layer then
in the uppermost one. Hence, the second statement of the Lema is proven.

Finally, the third statement follows on noting that the set $A$ described in
the second statement cannot be extended by any element to a set of 9 distinct
elements without a zerosum of length $\leq 3$.
\end{proof}

\begin{Lem}
\label{Lem:14points}
Let $A$ be a sequence of 14 points which does not contain a zerosum subset of
length $\leq 3$ or of length $\geq 12$. Then $A$ contains 7 distinct points,
each of which is taken twice. Moreover, there exists a multiset $A$ with these
properties, and it is unique up to linear equivalence.
\end{Lem}
\begin{proof}
The existence of $A$ is given by the following example.

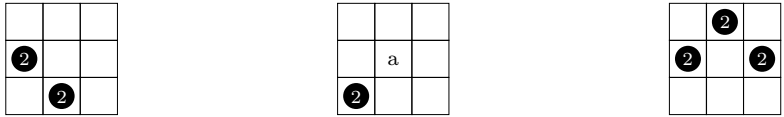
\begin{miniboard}
+ + +
2 + +
+ 2 +
\end{miniboard}\hspace*{\fill}
\begin{miniboard}
+ + +
+ a +
2 + +
\end{miniboard}\hspace*{\fill}
\begin{miniboard}
+ 2 +
2 + 2
+ + +
\end{miniboard}

It is easy to check that $A$ does not contain a zerosum of length $\leq
3$. Next, the sum of all elements in $A$ equals $(2, 2, 2)$; the inverse of
this element being the element maked $a$ in the picture above. This element is
neither 0 nor contained in $A$, hence, there is no zerosum of length $\geq
13$. Suppose that $a$ was the sum of two elements in $A$. Then $a$ is either
the sum of an element in the lowest layer with an element of the middle layer,
or the sum of two elements in the upper layer, but both possibilities are
easily dismissed. 

We now show that all sets of such 14 elements have indeed the form described
above. Thus, let $A$ be a set consisting of 14 elements, 8 of which are
distinct. Denote by $B$ the configuration described in 
Lemma~\ref{Lem:distinctfew}, part (ii). Then up to some linear transformation,
$A$ consists of 6 points of $B$ taken twice, and the 2 remaining points of $B$
taken once. Note that the sum of all elements in $B$ is 0, and the sum of any
two elements is non-zero, hence, the sum of all elements in $A$ is non-zero as
well. But $B$ is maximal among all sets of distinct elements without zerosum
subsets of length $\leq 3$, thus, every non-zero element, and in particular
the sum of all elements in $A$, can be represented as
the sum of 1 or 2 elements in $B$. Hence, by deleting 1 or 2 elements of $A$
we obtain a zerosum subset consisting of 12 or 13 elements.
\end{proof}

\begin{prop}
We have $D_2^*(\Z_3^3)=10$ 
\end{prop}
\begin{proof}
Consider the example 

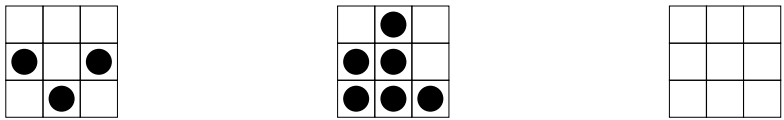
\begin{miniboard}
+ + +
X + X
+ X +
\end{miniboard}\hspace*{\fill}
\begin{miniboard}
+ X +
X X +
X X X
\end{miniboard}\hspace*{\fill}
\begin{miniboard}
+ + +
+ + +
+ + +
\end{miniboard}

Clearly there is no zerosum subset in the lowermost layer. Hence, if there are
two disjoint zero sum subsets, each of them must contain precisely 3 elements
of the second layer. Now consider the sum of all elements not contained in the
two zerosum subsets. This sum is equal to the sum of all elements in $A$, and
therefore $(2, 0, 0)$, on the other hand, it equals a subset sum of the
elements in the lowermost layer. However, $(2, 0, 0)$ cannot be represented by
elements in the lowermost layer. Hence, $D^*_2(\Z_3^3)\geq 10$. On the other
hand, $D^{3*}(\Z_3^3)=9$, thus among 10 points there is a zerosum subset of
length $\leq 3$, and among the remaining 7 points, there is always another
zerosum subset.
\end{proof}

\begin{prop}
We have $D^{*4}(\Z_3^3)=8$.
\end{prop}
\begin{proof}
Let $A$ be a set of 8 distinct elements which does not contain a zerosum
subset of length 4. Then $A$ cannot be contained in one plane, hence, we may
choose a basis of $\Z_3^3$ in $A$, which without loss is the standard
basis. Moreover, there has to be a sum in $A$, which we may assume,without 
any loss, to be $(0, 1, 1)$. Finally, we may change the third element of the
basis in such a way that the middle layer contains at least as many elements
as the upper layer, thus, we obtain the following picture, where at least two
elements of $A$ in the middle layer are not yet drawn.

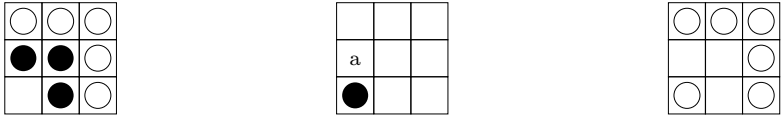
\begin{miniboard}
O O O
X X O
+ X O
\end{miniboard}\hspace*{\fill}
\begin{miniboard}
+ + +
a + +
X + +
\end{miniboard}\hspace*{\fill}
\begin{miniboard}
O O O
+ + O
O + O
\end{miniboard}

Suppose that $a=(1, 0, 1)$ is contained in $A$. Then several other elements
of $G$ can be excluded, since they would immediatelly give zerosum subsets of
size $\leq 4$, and we obtain the following situation.

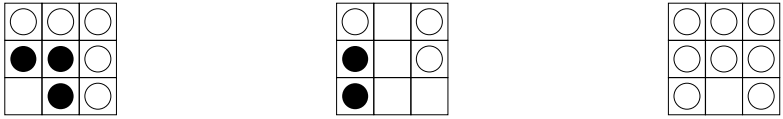
\begin{miniboard}
O O O
X X O
+ X O
\end{miniboard}\hspace*{\fill}
\begin{miniboard}
O + O
X + O
X + +
\end{miniboard}\hspace*{\fill}
\begin{miniboard}
O O O
O O O
O + O
\end{miniboard}

Hence, there is at most one element in the uppermost layer, that is, there are
at least two more in the middle layer. But any two elements in the middle row
of the middle layer would give a contradiction, hence, $(1, 2, 0)$ is in
$A$. But then no other element of $A$ could be in the middle layer, giving a
contradiction. Hence, $(1, 0, 1)$, and by symmetry $(1, 1, 0)$ are not
contained in $A$. Now assume that $(1, 0, 2)$ is in $A$. Then we obtain the
following situation.

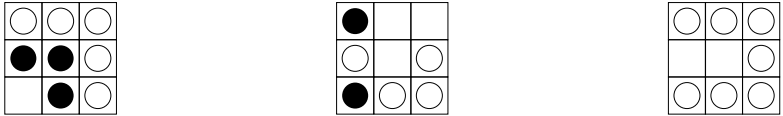
\begin{miniboard}
O O O
X X O
+ X O
\end{miniboard}\hspace*{\fill}
\begin{miniboard}
X + +
O + O
X O O
\end{miniboard}\hspace*{\fill}
\begin{miniboard}
O O O
+ + O
O O O
\end{miniboard}

By direct inspection we see that if both possible elements in the uppermost
layer are contained in $A$, then no additional element in the second layer
could be chosen, and $A$ would have at most 7 elements. Hence, two more
elements of $A$ in the middle layer are not yet shown. If $(1, 2, 2)$ was in
$A$, this is impossible, thus, we find that $(1, 1, 1)$ and $(1, 1, 2)$ are
both in $A$. Then we reach the following situation, which immediatelly implies
$|A|=7$, thus showing that $(1, 0, 2)$ and therefore by symmetry $(1, 2, 0)$,
are not in $A$ and hence, we obtain the following situation.

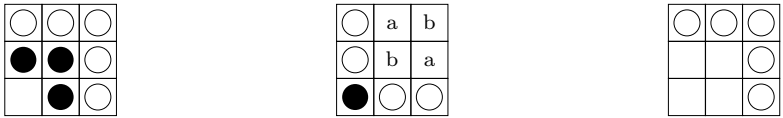
\begin{miniboard}
O O O
X X O
+ X O
\end{miniboard}\hspace*{\fill}
\begin{miniboard}
O a b
O b a
X O O
\end{miniboard}\hspace*{\fill}
\begin{miniboard}
O O O
+ + O
+ + O
\end{miniboard}

By assumption, there are two more elements in the middle layer, but on the
cells marked a and b, repsectively, there can be atmost one element of
$A$. Hence, without loss, we may assume that $(1, 1, 2)$ is in $A$, and, since
$(1, 0, 0)+(1, 1, 1)+(1, 1, 2)+(0, 1, 0)=(0, 0, 0)$, that $(1, 2, 2)$ is in
$A$ as well, that is, we reach the following situation.

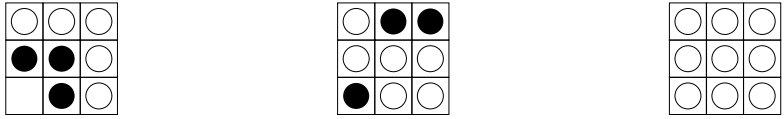
\begin{miniboard}
O O O
X X O
+ X O
\end{miniboard}\hspace*{\fill}
\begin{miniboard}
O X X
O O O
X O O
\end{miniboard}\hspace*{\fill}
\begin{miniboard}
O O O
O O O
O O O
\end{miniboard}

Clearly, $|A|=6$, contrary to our assumption, and we see that the initial
assumption on the existence of $A$ is wrong.
\end{proof}

\begin{prop}
\label{prop:D^4}
We have $D^4(\Z^3_3)=10$.
\end{prop}
\begin{proof}
The fact that $D^4(\Z^3_3)\geq 10$ is proven by the following example.

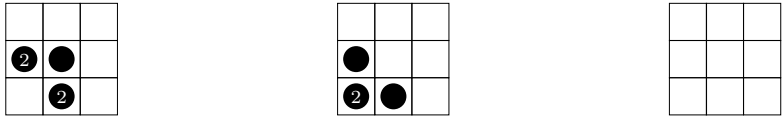
\begin{miniboard}
+ + +
2 X +
+ 2 +
\end{miniboard}\hspace*{\fill}
\begin{miniboard}
+ + +
X + +
2 X +
\end{miniboard}\hspace*{\fill}
\begin{miniboard}
+ + +
+ + +
+ + +
\end{miniboard}

Hence, it remains to show that every set of 10 elements contains a zerosum of
length at most 4. By means of contraiction, let $A$ be a set consisting of 10
elements of $G$ without a zerosum subset of 
size at most 4. Since $D^{*4}(\Z_3^3)=8$, at least 3 elements of $A$ are
repeated, and these 3 elements form a basis of $G$. Hence, we have the
following situation.

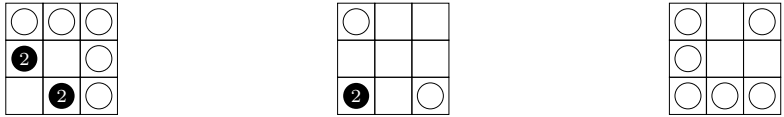
\begin{miniboard}
O O O
2 + O
+ 2 O
\end{miniboard}\hspace*{\fill}
\begin{miniboard}
O + +
+ + +
2 + O
\end{miniboard}\hspace*{\fill}
\begin{miniboard}
O + O
O + +
O O O
\end{miniboard}

Suppose that 2 of the three elements $(0, 1, 1)$, $(1, 0, 1)$ and $(1, 1, 0)$
are in $A$. Then we may suppose without loss that these elements are $(1, 1,
0)$ and $(1, 0, 1)$. Then we obtain the following.

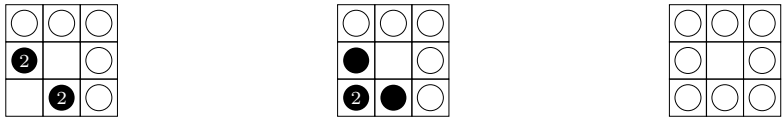
\begin{miniboard}
O O O
2 + O
+ 2 O
\end{miniboard}\hspace*{\fill}
\begin{miniboard}
O O O
X + O
2 X O
\end{miniboard}\hspace*{\fill}
\begin{miniboard}
O O O
O + O
O O O
\end{miniboard}

Since $2\cdot(1, 1, 0)+(1, 0, 0) + (0, 1, 0)=(0, 0, 0)$, both $(1, 1, 0)$ and
$(1, 0, 1)$ are taken at most once. Hence, there are two elements in $A$ of
the form $(a, 1, 1)$. If they are distinct, we have $(a, 1, 1)+(b, 1, 1)+(x,
1, 0) + (y, 0, 1)= (0, 0, 0)$ for appropriate valeus $x, y\in\{0, 1\}$, hence,
there is one value which is taken twice. However, all three remaining choices
lead to contradictions, and we conclude that at most one of $(0, 1, 1)$, $(1,
0, 1)$ and $(1, 1, 0)$ is contained in $A$; without loss we may assume that
$(1, 0, 1)$ and $(1, 1, 0)$ are not in $A$. Next we note that any two of $(1,
1, 2)$, $(1, 2, 1)$ and $(2, 1, 1)$ together with the elements already placed
give a zerosum subset of size 3, hence, at most one of these elements can be
contained in $A$. By symmetry we may suppose that $(1, 2, 1)$ is not in $A$;
moreover, if $(0, 1, 1)$ was not in $A$, then we may also assume that $(1, 1,
2)$ is not in $A$.

We shall now assume that $(0, 1, 1)$ is not in $A$. Then we have the following
situation.

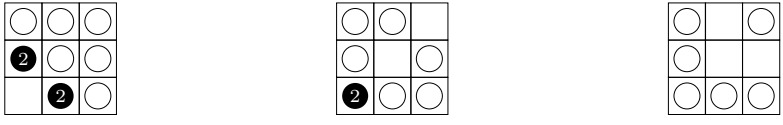
\begin{miniboard}
O O O
2 O O
+ 2 O
\end{miniboard}\hspace*{\fill}
\begin{miniboard}
O O +
O + O
2 O O
\end{miniboard}\hspace*{\fill}
\begin{miniboard}
O + O
O + +
O O O
\end{miniboard}

There are at most 8 elements in the lower two layers, hence, there are at
least two elements in the uppermost layer; in particular, $(1, 1, 1)$ is not
in $A$. Suppose that $(1, 2, 2)$ does not occur twice in $A$. Then there are
at least 3 elements of $A$ in the uppermost layer, and since $(2, 1, 2)$ and
$(2, 2, 1)$ cannot both occur in $A$, we deduce that $(2, 1, 1)$ and one of
$(2, 1, 2)$ and $(2, 2, 1)$ is in $A$; without loss we may assume the former,
and obtain the following situation.

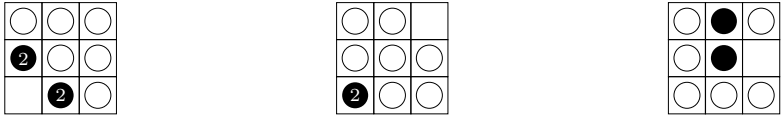
\begin{miniboard}
O O O
2 O O
+ 2 O
\end{miniboard}\hspace*{\fill}
\begin{miniboard}
O O +
O O O
2 O O
\end{miniboard}\hspace*{\fill}
\begin{miniboard}
O X O
O X +
O O O
\end{miniboard}

Since $(2, 1, 2)+(1, 2, 2)+2\cdot(0, 0, 1)=(0, 0, 0)$, we see that $(1, 2, 2)$
is not contained in $A$, and we find that the elements in the uppermost layer
are contained twice in $A$. But then we obtain the contradiction $2\cdot(2, 1,
2)+(2, 1, 1)+(0, 0, 1)=(0, 0, 0)$.

Hence, our initial assumption that $(0, 1, 1)$ was not in $A$ is false, and we
obtain the following situation.

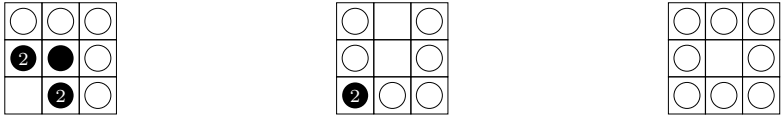
\begin{miniboard}
O O O
2 X O
+ 2 O
\end{miniboard}\hspace*{\fill}
\begin{miniboard}
O + O
O + O
2 O O
\end{miniboard}\hspace*{\fill}
\begin{miniboard}
O O O
O + O
O O O
\end{miniboard}

Note that $(0, 1, 1)$ cannot be repeated in $A$, because $D^{4}(\Z_3^2)=6$,
thus, there are 3 elements of $A$ on the remaining three cells. But $(1, 1,
1)$ and $(2, 1, 1)$ cannot be simultaneously in $A$, thus, $(1, 1, 2)$ is in
$A$. Because of $(1, 0, 0)+(0, 1, 0)+(1, 1, 1)+(1, 1, 2)=(0, 0, 0)$, $(1, 1,
1)$ cannot be contained in $A$, and therefore $(2, 1, 1)$ has to be contained
in $A$. But then we obtain the zerosum $(2, 1, 1)+(1, 1, 2)+(0, 1, 0)$, and
obtain a contradiction proving our claim.
\end{proof}

\begin{prop}
\label{prop:D_2}
We have $D_2(\Z_3^3)=11$.
\end{prop}
\begin{proof}
Let $A$ be a subset of $\Z_3^3$ containing 11 elements. Then $A$ contains a
zerosum subset of size $\leq4$, and the complement of this zerosum is a set
with $\geq 7$ elements, which therefore contains another zerosum
subset. Hence, $D_2(\Z_3^3)\leq 11$. On the other hand, the inequality
$D_2(\Z_3^3)\geq11$ is proven by the following example.

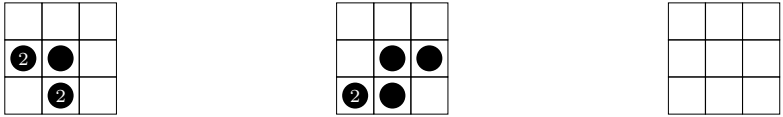
\begin{miniboard}
+ + +
2 X +
+ 2 +
\end{miniboard}\hspace*{\fill}
\begin{miniboard}
+ + +
+ X X
2 X +
\end{miniboard}\hspace*{\fill}
\begin{miniboard}
+ + +
+ + +
+ + +
\end{miniboard}

In fact, every zerosum subset contains either no elements of the middle layer,
or precisely 3, thus, if there were to distinct zerosum subsets, one of them
has to be contained in the lowermost layer. Obviously, it has to contain all
points of this layer. But there is no zerosum subset of size three in the
middle layer, and we find that this set does not contain two disjoint zerosum
subsets, which proves $D_2(\Z_3^3)\geq11$.
\end{proof}

\begin{prop}
\label{prop:D_3}
We have $D_3(\Z_3^3) = 15$.
\end{prop}
\begin{proof}
The upper bound $D_3(\Z_3^3)\leq 15$ follows from
Proposition~\ref{prop:D^4} and \ref{prop:D_2}, whereas the lower bound
follows from the configuration given in Lemma~\ref{Lem:14points}.
\end{proof}

\begin{prop}
\label{prop:D^5}
We have $D^5(\Z_3^3)=9$.
\end{prop}
\begin{proof}
The lower bound $D^5(\Z_3^3)\geq 8$ is given by the following example.

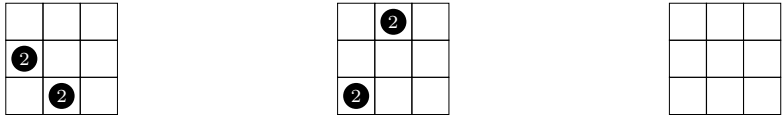
\begin{miniboard}
+ + +
2 + +
+ 2 +
\end{miniboard}\hspace*{\fill}
\begin{miniboard}
+ 2 +
+ + +
2 + +
\end{miniboard}\hspace*{\fill}
\begin{miniboard}
+ + +
+ + +
+ + +
\end{miniboard}

Let $A$ be a set consisting of 9 elements without a zerosum subset of length
5. Suppose first that there are at most 5 distinct elements in $A$. Then there
are at least 4 elements twice in $A$, and these elements cannost be in one
plane, hence, there is a basis of elements taken twice in $A$. We therefore
obtain the following.

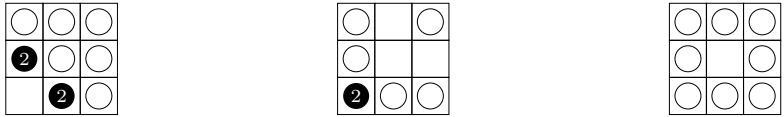
\begin{miniboard}
O O O
2 O O
+ 2 O
\end{miniboard}\hspace*{\fill}
\begin{miniboard}
O + O
O + +
2 O O
\end{miniboard}\hspace*{\fill}
\begin{miniboard}
O O O
O + O
O O O
\end{miniboard}

If $(2, 1, 1)\in A$, then there are no further elements in the middle layer,
which would imply $|A|\leq 8$; thus, by symmetry, we have $(1, 1, 2), (1, 2,
1), (2, 1, 1)\not\in A$. But then $A$ could only have $(1, 1, 1)$ as a
possible element left, and therefore $|A|\leq 8$. 

Hence, we may assume that $A$ contains 6 distinct elements, and therefore
there is some basis $\{a, b, c\}\in A$, such that $a+b\in A$, that is, we have
the following situation.

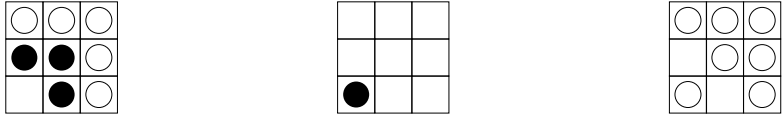
\begin{miniboard}
O O O
X X O
+ X O
\end{miniboard}\hspace*{\fill}
\begin{miniboard}
+ + +
+ + +
X + +
\end{miniboard}\hspace*{\fill}
\begin{miniboard}
O O O
+ O O
O + O
\end{miniboard}

Moreover, we can choose the third base element in such a way that it is either
twice in $A$, or that no element outside the lowest plane is twice in
$A$. Consider first the case that $A$ occurs twice. 

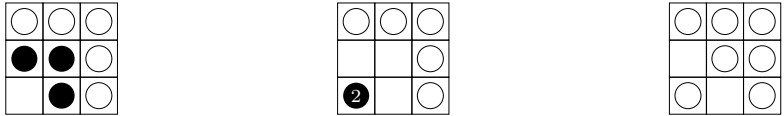
\begin{miniboard}
O O O
X X O
+ X O
\end{miniboard}\hspace*{\fill}
\begin{miniboard}
O O O
+ + O
2 + O
\end{miniboard}\hspace*{\fill}
\begin{miniboard}
O O O
+ O O
O + O
\end{miniboard}

Note that $(0, 1, 1)$ cannot be contained twice in $A$, and that $(0, 0, 1)$
and $(0, 1, 0)$ cannot be both twice in $A$, that is, there are at most 4
elements in the lowest layer. Moreover, all elements in the middle layer not
yet depicted can occur at most once, and not all three empty places can be
taken, thus, there is an eleemnt in the uppermost layer; without loss we may
suppose that this element is $(2, 0, 1)$. We then obtain the following.

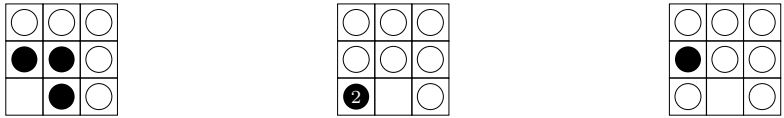
\begin{miniboard}
O O O
X X O
+ X O
\end{miniboard}\hspace*{\fill}
\begin{miniboard}
O O O
O O O
2 + O
\end{miniboard}\hspace*{\fill}
\begin{miniboard}
O O O
X O O
O + O
\end{miniboard}

Here, $(2, 0, 1)$ can only be once in $A$, since otherwise we had the zerosum
$2\cdot(2, 0, 1)+2\cdot(1, 0, 0) + (0, 0, 1)$, thus, there are at most 2
points in the uppermost layer. Hence, $(1, 1, 0)\in A$, which implies that
$(2, 1, 0)\not\in A$, and we find that $|A|\leq 8$. Thus, from now on, we
shall assume that all elements outside the lowest layer occur only once in $A$.

Suppose that $(1, 0, 1)\in A$. Then we have the following.

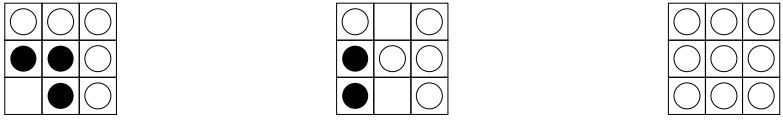
\begin{miniboard}
O O O
X X O
+ X O
\end{miniboard}\hspace*{\fill}
\begin{miniboard}
O + O
X O O
X + O
\end{miniboard}\hspace*{\fill}
\begin{miniboard}
O O O
O O O
O O O
\end{miniboard}

Obviously, there can only be 4 elements outside the first layer, and at most 4
inside the first layer, which gives a contradiction to $|A|=9$. 

Next, suppose that $(1, 1, 2)\in A$. Then we obtain

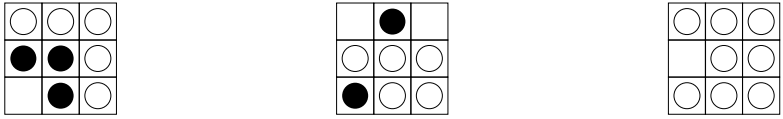
\begin{miniboard}
O O O
X X O
+ X O
\end{miniboard}\hspace*{\fill}
\begin{miniboard}
+ X +
O O O
X O O
\end{miniboard}\hspace*{\fill}
\begin{miniboard}
O O O
+ O O
O O O
\end{miniboard}

Again, there are only 5 places for elements outside the first layer, and $(1,
0, 2)$ and $(1, 2, 2)$ cannot occur at the same time, thus, $(1, 1, 2)\not\in
A$ as well.

Hence, we are led to the following constellation.

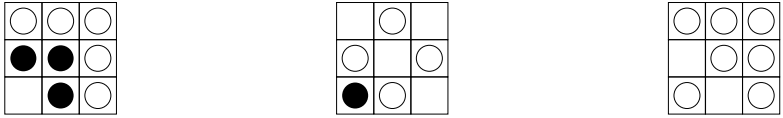
\begin{miniboard}
O O O
X X O
+ X O
\end{miniboard}\hspace*{\fill}
\begin{miniboard}
+ O +
O + O
X O +
\end{miniboard}\hspace*{\fill}
\begin{miniboard}
O O O
+ O O
O + O
\end{miniboard}

There are only 7 places left for elements outside the first layer. Moreover,
among these $(1, 1, 1)$ and $(1, 2, 2)$ as well as $(1, 2, 0)$ and $(1, 0, 2)$
are mutually exclusive, and we conclude that $(2, 1, 0), (2, 0, 1)\in A$, and,
without loss, $(1, 2, 0)\in A$. But then we obtain the zerosum $(2, 0, 1)+(1,
2, 0)+(0, 1, 1)+(0, 0, 1)$, and this contradiction proves our claim.
\end{proof}

\begin{prop}
\label{Prop:3k+5Char}
For $k\geq 3$ we have $D_k(\Z_3^3)=3k+6$. Moreover, the set of all multisets
$A$ of size $3k+5$ which do not have $k$ disjoint zerosum subsets can be
constructed as follows: Take all sets $B=\{b_1, \ldots, b_7\}$ of 7 distinct
points without a subsum of length $\leq 3$, such that the multiset $C$
obtained from $B$ by taking each point twice does not contain a zerosum subset
of length $\geq 12$. Choose a partition $k-3=\kappa_1+\ldots+\kappa_7$, and
set $A=C\cup\{b_1^{3\kappa_1}, \ldots, b_7^{3\kappa_7}\}$.
\end{prop}
\begin{proof}
We first show that none of the sets described here contain $k$ disjoint zersosum
subsets. In fact, since there is no zersum of length $\leq 3$ in $B$, every
zerosum of length 3 in $A$ must contain the same element three times, thus,
any collection of disjoint zerosum subsets can contain at most $k-2$ zerosums
of length 3. next, note that the sum of all elements of $A$ equals twice the
sum of all elements of $B$, and that the set of elements representable by 0, 1
or 2 elements of $A$ is equal to the set of elements representable by 0, 1, or
2 elements of $B\cup B$; thus, as in the proof of Lemma\ref{Lem:14points} we see
that there is no zerosum of length $\geq 3k+3$. Hence, every collection of
disjoint zerosum subsets can contain $3k+2$ points at most. Thus, such a
collection contains no zerosum subset of length $\leq 2$, at most $k-3$ of
length 3, and alltogether consists of at most $3k+2$ points, which implies
that the total number of zerosum subsets is $\leq k-1$.

Now we show that there are no examples different from the one described here.
Let $A$ be a multiset consisting of $3k+5$ elements of $\Z_3^3$ which does not
have $k$ disjoint zerosum subsets. If there is a
zerosum of length $\leq 2$ in $A$, then removing this zerosum yields $k-1$
zerosums in the remaining $3(k-1)+6$ points, which gives a contradiction. Next
suppose there is an element repeated 4 times in $A$. Then we can remove this
element 3 times, and may assume by induction that the new set has the form
described. hence, it suffices to consider the case that every element in $A$
occurs at most 3 times.

Suppose there is one element $a$ occurring once. Then
we remove as many zerosum subsets of length 3 from $A$ as possible without
removing this point. We end up with a set $B$ of 14 or 17 points, which
contains one element precisely once, and does not contain 3 resp. 4 disjoint
zerosum subsets. However, we already saw that this is impossibel for a set of
14elements. If $|B|=17$, then either $B$ contains one element three times,
contradicting the assumption that we removed all zerosums of length three
avoiding $a$, or there are at least 9 distinct points in $B$. In the latter
case let $\ell\leq 8$ be the number of points occurring twice in $B$. Collect
each point which is twice in $B$ and as many points that occur once in $B$
necessary to reach 9 points in a set $C$. Then $C$ contains a zerosum of
length 3, we claim that removing this zerosum of $B$ yields a set which
contains one element only once. This is clear if $b<8$, for then some element
occurring once in $B$ is not contained in $C$, and can therefore also not be a
part of the zerosum removed. If on the other hand $b=8$, then at least two
points occurring twice in $B$ are removed once, thusthere are at least 2
points in the new set with multiplicity 1. In any case, we have removed $k-3$
zerosum subsets of size 3 from the beginning set $A$, and ended up in a set of
14 elements, one of which occurring only once. Hence, there are 3 more
zerosums, contradicting the assumption that $A$ does not contain $k$ zerosum
subsets. 

Next, suppose there is a zerosum of length 3 which consists of distinct
elements. Then we can remove this zerosum once or twice to reach a situation
with one element occurring precisely once, a situation we just dealt with. In
particular, there are at most 8 distinct elements in $A$. If there are only 7
distinct elements, we reached the position described in
Lemma~\ref{Lem:14points} 
and are done. Otherwise we have a set with 17 elements, 8 distinct ones among
them forming the set described in Lemma~\ref{Lem:distinctfew} (ii), and
one point occurs 
three times, whereas the other points occur exactly twice. Removing the three
times repeated point once, we obtain a set of 16 elements with sum 0, thus,
among the 15 remaining points there are 3 disjoint zerosums, whereas the
complement of these three sums constitue a fourth one, yielding $k$ disjoint
zerosums for the original set $A$.

Hence, our claim follows.
\end{proof}

\section{Proof of Theorem~\ref{thm:D=MZ3}}

\begin{Lem}
In every set of 5 distinct elements of $\Z_3^3$ there is either a
zerosum of length $\leq 3$, or there are 3 elements $x, y, z$ satisfying
the equation $x+y=z$.
\end{Lem}
\begin{proof}
It is easy to check that the analogous statement in two dimensions
holds true for all sets of 3 elements, hence, we may assume that $A$
does not contain 3 points in any plane passing through the origin, and
we obtain the following situation.

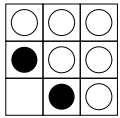
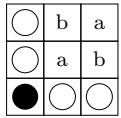
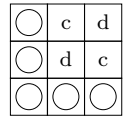
\begin{miniboard}
O O O
X O O
+ X O
\end{miniboard}\hspace*{\fill}
\begin{miniboard}
O b a
O a b
X O O
\end{miniboard}\hspace*{\fill}
\begin{miniboard}
O c d
O d c
O O O
\end{miniboard}

Without loss we may suppose that the middle plane contains at least
one other element of $A$.
If one of the cells marked a is in $A$, then none of the cells marked
b is in $A$, and vice versa. On the other hand, if both cells marked a
are in $A$, we would obtain a zerosum of length 3, and similar for b,
thus, precisely one cell in the middle plane is contained in
$A$. Moreover, if a cell marked a is in $A$, the cells marked d cannot
be in $A$, and similarly with b, and we are left with the
possibilities that there is precisely one element a and one c, or one
b and one d.

If $(1, 1, 1)\in A$, both elements marked c yield zerosums of length 3,
similarly, if $(2, 1, 1)\in A$. If $(1, 1, 2)\in A$, then the only
remaining possibility for the last element is $(2, 2, 2)$, but $(1, 1,
2)+(2, 2, 2)=(0, 0, 1)\in A$, and the same argument applies to the case
that $(2, 1, 1)$ is in $A$. Hence, no such set $A$ can exist.
\end{proof}

\begin{Theo}
\label{thm:Af}
Let $n$ be an integer coprime to 6. Then there does not exist a multiset
$A\subseteq\Z_3^3$ with 10 elements together with a function $f:A\to
\Z_n$, such that $A$ does not contain 2 disjoint zerosum subsets, that for
every zerosum subset $B$ of $A$ we have $\sum_{b\in B} f(b) = 1$, and that
there is some element $a\in A$ with $3f(a)=1$.
\end{Theo}
\begin{proof}[Proof of the theorem]
Suppose that $A$ and $f$ are as in the statement. Observe that $A$
does not have a zerosum of length  $\leq 3$ or $\geq 8$. In
particular, there are at most 8 distinct elements in $A$, that is,
there are at least 3 elements occurring twice. We distinguish
cases according to the constellation of these elements.

(i) Suppose that there exist 3 elements $a, b, c$ occurring twice within one
plane passing through the origin. Without loss we may suppose that
$c=a+b$. Then we have the zerosums $2a+2b+c=a+b+2c=0$, thus
$f(a)+f(b)=f(c)$ and $3f(c)=1$, in particular, $f(c)\neq 0$. Suppose that
in $A$ there are elements $x_i$ outside this 
plane adding up to $2a$. Then we have the zerosums $a+\sum x_i =
2a+b+2c+\sum x_i=0$, which implies the equation
\[
f(a) + \sum f(x_i) = 4f(a) + 3f(b) + \sum f(x_i),
\]
which in turn implies $f(c)= 0$, a contradiction. Next suppose that
$2a+2b$ can be represented as a sum of elements $x_i$ of $A$ outside the
given plane. Then we have the zerosums $2a+2b+2c+\sum x_i$ and $c+\sum x_i$,
which in the same way implies $f(c)=0$. If $2a+b$ can be represented,
we have the zerosums $a+2b+\sum x_i$ and $2a+2c+\sum x_i$, which implies
$3f(a)=0$, that is, $f(a)=0$.
If $a$ can be represented, we obtain the zerosums $2a+\sum x_i$ and
$b+2c+\sum x_i$, which implies $f(b)=0$. Noting that not both $f(a)$ and
$f(b)$ can vanish, since otherewise $f(c)$ would be zero, we can
summarize these considerations as follows: None of $2a, 2b, 2c$ can be
represented as a sum of elements of $A$ outside the plane spanned by
$a$ and $b$, and if one of $a, a+2b$ is represented in such a way, none
of $2a+b, b$ is, and vice versa.

(ii) As an immediate consequence we obtain that $A$ cannot consist of
5 elements, each taken twice. In fact, without loss, we have the
following situation.

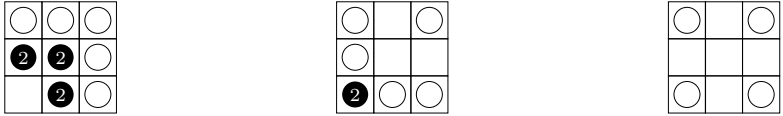
\begin{miniboard}
O O O
2 2 O
+ 2 O
\end{miniboard}\hspace*{\fill}
\begin{miniboard}
O + O
O + +
2 O O
\end{miniboard}\hspace*{\fill}
\begin{miniboard}
O + O
+ + +
O + O
\end{miniboard}

If there is an element $x$ in the middle layer occurring twice in $A$,
$2\cdot (1, 0, 0) + x$ and $(1, 0, 0)+2x$ are elements in the lowest
layer; if there is an element $x$ in the upermost layer occurring
twice, both $(1, 0, 0)+x$ and $2\cdot(1, 0, 0) + 2x$ are elements in the
lowest layer. In any case, there is some $y$ in the lowest layer, such
that both $y$ and $2y$ can be represented as a sum of elements outside
the lowermost layer. If $y=0$, we have a zerosum of length 2, if
$y=(0, 1, 2)$, both $(0, 1, 2)$ and $(0, 2, 1)$ are representable,
which gives a contradiction, and in all other cases one of $(0, 0,
2)$, $(0, 2, 0)$ and $(0, 2, 2)$ is representable, which gives also a
contradiction. 

(iii) Suppose that $(0, 0, 1), (0, 1, 0)$ and $(0, 1, 1)$ occur twice,
and that none of the equations $x+(0, 1, 0)=y, x+(0, 0, 1)=y$ is
solvable with $x, y\in A$. Without loss we may assume that $(1, 0, 0)$
is in $A$, and that, if some other element occurs twice in $A$, then
so does $(1, 0, 0)$. We obtain the following situation.

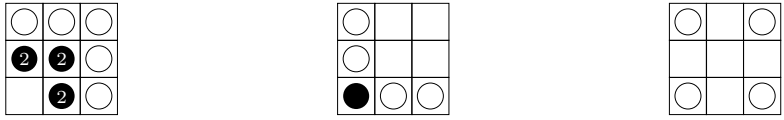
\begin{miniboard}
O O O
2 2 O
+ 2 O
\end{miniboard}\hspace*{\fill}
\begin{miniboard}
O + +
O + +
X O O
\end{miniboard}\hspace*{\fill}
\begin{miniboard}
O + O
+ + +
O + O
\end{miniboard}

If there are two elements in the uppermost layer, we may assume
without loss that $(2, 1, 0)\in A$. Then $(2, 0, 1)$ and $(2, 2, 1)$ 
cannot be in $A$ since otherwise $(0, 1, 0)$ and one of $(0, 0, 1)$,
$(0, 2, 1)$ is representable by elements outside the lowest layer; and
$(2, 1, 1)$ and $(2, 1, 2)$ cannot be contained in $A$, since we would
obtain a solution of the equation $x+(0, 1, 0=y$ in $A$. Hence, there
is at most one element of $A$ in the uppermost layer. On the other
hand, there can be at most 3 elements in the middle layer, thus, there
is precisely one element in the uppermost layer, and 3 elements in the
middle layer, two of which are $(1, 0, 0)$. Moreover, $(2, 1, 1)\not\in
A$, since otherwise there would be no place left for the last element
in the middle layer, and $(1, 2, 2)\not\in A$, since otherwise $(0, 2,
2)$ was representable. We therefore get the following.

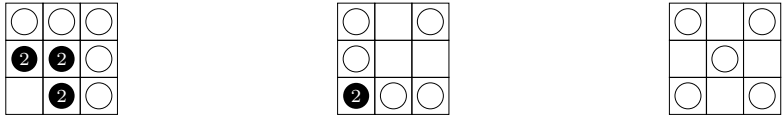
\begin{miniboard}
O O O
2 2 O
+ 2 O
\end{miniboard}\hspace*{\fill}
\begin{miniboard}
O + O
O + +
2 O O
\end{miniboard}\hspace*{\fill}
\begin{miniboard}
O + O
+ O +
O + O
\end{miniboard}

We now check that $(2, 1, 2)\in A$ is impossible, since there is no
place for the last element in the middle layer, and we may assume
without loss that $(2, 1, 0)$ is the unique element in the uppermost
layer.

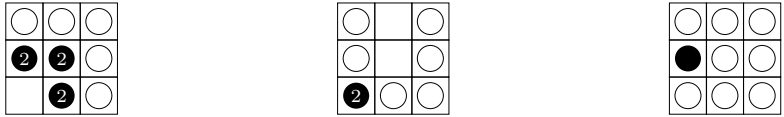
\begin{miniboard}
O O O
2 2 O
+ 2 O
\end{miniboard}\hspace*{\fill}
\begin{miniboard}
O + O
O + O
2 O O
\end{miniboard}\hspace*{\fill}
\begin{miniboard}
O O O
X O O
O O O
\end{miniboard}

The sum of all elements in $A$ is therefore $(2, 0, 1)$ or $(2, 0,
2)$, and both these elements can be represented as a sum of 1 or 2
elements of $A$; thus, $A$ contains a zerosum of length 8 or 9, which
gives a contradiction.

(iv) Now suppose that $(0, 0, 1), (0, 1, 0)$ and $(0, 1, 1)$ occur
twice. By the preceding argument we know that for some elements $x, y\in
A$ outside the lowermost plane we have $x+(0, 1, 0)=y$ or $x+(0, 0,
1)=y$, and we may therefore assume without loss that both $(1, 0, 0)$ and
$(1, 0, 1)$ are contained in $A$. We obtain the following situation.

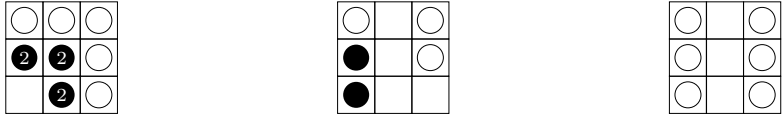
\begin{miniboard}
O O O
2 2 O
+ 2 O
\end{miniboard}\hspace*{\fill}
\begin{miniboard}
O + O
X + O
X + +
\end{miniboard}\hspace*{\fill}
\begin{miniboard}
O + O
O + O
O + O
\end{miniboard}

Suppose first that $(1, 0, 0)$ is not twice in $A$. Then $(1, 2, 0)$
is in $A$, since $(1, 0, 1)$ cannot be twice in $A$, and therefore all
remaining elements are in the plane $(s, 1, t)$, which would imply
that the sum of all elements outside the lowermost plane is of the form
$(0, 2, t)$, which for $t\neq 1$ yields a contradiction, whereas for
$t=1$ we deduce that the sum of all elements in $A$ is $(0, 0, 2)$,
which is a sum of two elements in $A$, thus, $A$ contains a zerosum of
length 8, which also yields a contradiction. Hence, we conclude that
$(1, 2, 0)$ is in $A$, and obtain the following.

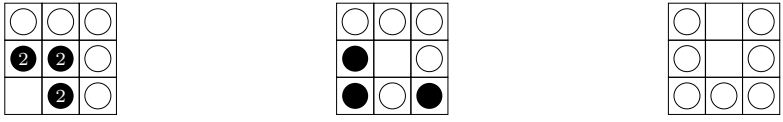
\begin{miniboard}
O O O
2 2 O
+ 2 O
\end{miniboard}\hspace*{\fill}
\begin{miniboard}
O O O
X + O
X O X
\end{miniboard}\hspace*{\fill}
\begin{miniboard}
O + O
O + O
O O O
\end{miniboard}

Now $(0, 2, 1)$ is a sum of elements outside the lowermost plane,
hence, $(0, 1, 0)$ and $(0, 1, 2)$ are not, which implies that $(1, 1,
1), (2, 1, 1), (2, 1, 2),\not\in A$, and that $(1, 2, 0)$ is only once
in $A$, thus, $A$ consists only of 9 elements. This contradiction
shows that $(1, 0, 0)$ occurs twice in $A$, and we find the following.

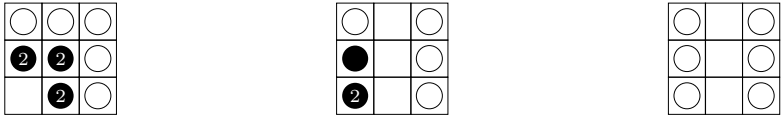
\begin{miniboard}
O O O
2 2 O
+ 2 O
\end{miniboard}\hspace*{\fill}
\begin{miniboard}
O + O
X + O
2 + O
\end{miniboard}\hspace*{\fill}
\begin{miniboard}
O + O
O + O
O + O
\end{miniboard}

We now can represent $(0, 0, 1)$ as a sum of elements outside the
lowermost layer, hence, $(0, 1, 0)$ and $(0, 1, 2)$ cannot be
represented this way. The remaining element of $A$ has the form $(s,
1, t)$. If $t\neq 1$, we add $(1, 0, 0)$ once or twice, to obtain a
forbidden sum in the lowermost layer, whereas if $t=1$, we add $(1, 0,
1)$, and, if necessary, $(1, 0, 0)$ to obtain a forbidden sum as well.
Hence, in any case, we obtain a contradiction.

Thus, we see that the equation $x+y=z$ is not solvable among distinct
elements occurring twice in $A$, in particular, we may assume that
$(1, 0, 0)$, $(0, 1, 0)$ and $(0, 0, 1)$ all occur twice in $A$.

(v) Suppose that there are 4 elements occurring twice in $A$. We
already saw that three of them may be taken to be $(0, 0, 1), (0, 1,
0)$ and $(1, 0, 0)$, and that the fourth cannot be contained in one of
the planes generated by this three. Hence, there are 8 possibilities
left. These 8 points fall into 4 equivalence classes under rotation
around the spatial diagonal, and taking the fourth point to be
$(2, 2, 2)$ would yield a zerosum of length 8. Moreover, a linear
transformation fixing the plane $(0, s, t)$ and interchanging $(1, 0,
0)$ and $(1, 1, 1)$ shows that the case of the fourth point being $(1,
1, 1)$ is equivalent to the case $(1, 2, 2)$. Hence, up to symmetry we
may suppose that the fourth point is $(1, 1, 1)$ or $(1, 1, 2)$.

(vi) Suppose that $(0, 0, 1), (0, 1, 0)$, $(1, 0, 0)$ and $(1, 1, 1)$
occur twice in $A$. Then we have the following.

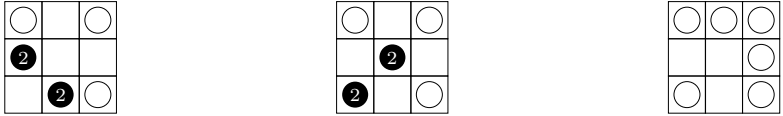
\begin{miniboard}
O + O
2 + +
+ 2 O
\end{miniboard}\hspace*{\fill}
\begin{miniboard}
O + O
+ 2 +
2 + O
\end{miniboard}\hspace*{\fill}
\begin{miniboard}
O O O
+ + O
O + O
\end{miniboard}

Up to symmetry all elements of $\Z_3^3$ except $(0, 1, 2)$, $(0, 2, 2)$
and $(1, 2, 2)$ can be represented as the sum of at most 2 elements already
depicted, thus we deduce that the sum of all elements of $A$ is equal
to one of these, for otherwise we would obtain a zerosum of length $\geq
8$. Substracting the elements depicted we find that the sum of the two
remaining elements equals $(2, 0, 1)$, $(2, 1, 1)$ or $(0, 1, 1)$.

Suppose first the sum is $(2, 0, 1)$. Deleting all elements $x$ such
that $(2, 0, 1)-x$ is impossible, only the following two possilities remain.

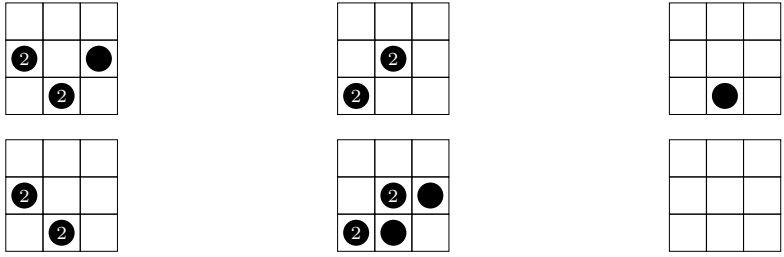
\begin{miniboard}
+ + +
2 + X
+ 2 +
\end{miniboard}\hspace*{\fill}
\begin{miniboard}
+ + +
+ 2 +
2 + +
\end{miniboard}\hspace*{\fill}
\begin{miniboard}
+ + +
+ + +
+ X +
\end{miniboard}

\begin{miniboard}
+ + +
2 + +
+ 2 +
\end{miniboard}\hspace*{\fill}
\begin{miniboard}
+ + +
+ 2 X
2 X +
\end{miniboard}\hspace*{\fill}
\begin{miniboard}
+ + +
+ + +
+ + +
\end{miniboard}

The same argument applied to the second case yields only one case.

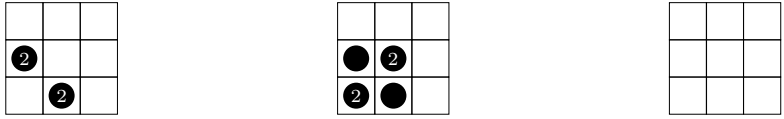
\begin{miniboard}
+ + +
2 + +
+ 2 +
\end{miniboard}\hspace*{\fill}
\begin{miniboard}
+ + +
X 2 +
2 X +
\end{miniboard}\hspace*{\fill}
\begin{miniboard}
+ + +
+ + +
+ + +
\end{miniboard}

Finally, the last sum gives two cases, which are symmetric to each
other, thus, we only have to consider the following. 

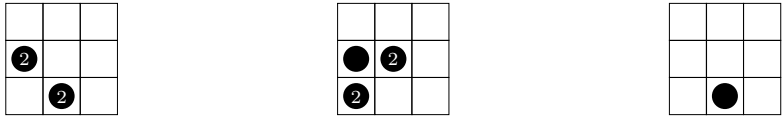
\begin{miniboard}
+ + +
2 + +
+ 2 +
\end{miniboard}\hspace*{\fill}
\begin{miniboard}
+ + +
X 2 +
2 + +
\end{miniboard}\hspace*{\fill}
\begin{miniboard}
+ + +
+ + +
+ X +
\end{miniboard}

Here, the sum of
all elements different from $(1, 1, 1)$ is 0, thus, there is a zerosum
of length 8, and it suffices to consider the 3 previous cases.

We set $x=(0, 0, 1)$, $y=(0, 1, 0)$, $z=(1, 0, 0)$, $w=(1, 1,
1)$. Then we have the zerosums $2x+2y+2z+w, x+y+z+2w$, which imply
that $f(w)=f(x)+f(y)+f(z)$ and $3f(w)=1$. Moreover, if $a=(r, s, t)$ is
one of the other elements of $A$, we have the zerosum
$a+[-r]z+[-s]y+[-t]x$, similarly, we have the zerosums
$a+w+[-1-r]z+[-1-s]y+[-1-t]x$ and
$a+2w+[-2-r]z+[-2-s]y+[-2-t]x$, and we obtain the equations
\begin{eqnarray*}
(1-[-r]+[-1-r]) f(z) + (1-[-s]+[-1-s]) f(y) + (1-[-t]+[-1-t]) f(x) & =
& 0\\
(2-[-r]+[-2-r]) f(z) + (2-[-s]+[-2-s]) f(y) + (2-[-t]+[-2-t]) f(x) & =
& 0\\
\end{eqnarray*}
Every coefficient is divisible by 3, and in the interval $[0, 5]$,
hence, dividing by 3 we obtain equations with all coefficients 0 and
1. We can apply this argument to both points in $A$ occurring only
once as well as to their sum, thus, in each of the four cases we
obtain a $3\oplus 6$-matrix with the 
property that $(f(z), f(y), f(x))$ is in the kernel of this
matrix. To compute these matrices, note that
\[
\frac{1-[-a]+[-1-a]}{3} = \begin{cases}0, &a=1,2\\1,&a=0\end{cases},
\quad
\frac{2-[-a]+[-2-a]}{3} = \begin{cases} 0, & a=1\\ 1, & a=0, 2\end{cases}
\]
for all $a\in\Z_3$. We therefore obtain the matrices
\[
\begin{pmatrix}
1 & 0 & 0\\ 1 & 1 & 0\\ 0 & 0 & 1\\ 1 & 0 & 1\\ 0 & 1 & 0\\ 1 & 1 & 0
\end{pmatrix}, \qquad
\begin{pmatrix}
0 & 0 & 1\\ 0 & 0 & 1\\ 0 & 0 & 0 \\ 0 & 1 & 0\\ 0 & 1 & 0\\ 1 & 1 & 0
\end{pmatrix}, \qquad
\begin{pmatrix}
0 & 0 & 1\\ 0 & 0 & 1\\ 0 & 1 & 0\\ 0 & 1 & 0\\ 0 & 0 & 0\\ 1 & 0 & 0
\end{pmatrix}
\]
Obviously, all these matrices have rank 3, which implies that
$f(x)=f(y)=f(z)=0$, and we obtain the contradiction
\[
1=3f(w)=3(f(x)+f(y)+f(z)) = 0.
\]

(vii) Suppose that $(0, 0, 1)$, $(0, 1, 0)$, $(1, 0, 0)$ and $(1, 1,
2)$ occur twice in $A$. Then we have the following situation.

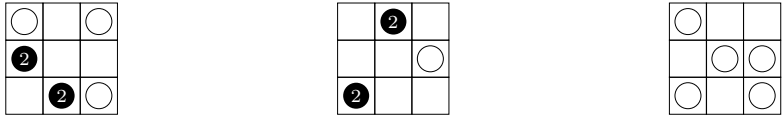
\begin{miniboard}
O + O
2 + +
+ 2 O
\end{miniboard}\hspace*{\fill}
\begin{miniboard}
+ 2 +
+ + O
2 + +
\end{miniboard}\hspace*{\fill}
\begin{miniboard}
O + +
+ O O
O + O
\end{miniboard}

As in the previous argument we find that $(r, s, t)\in A$ implies the
equations
\begin{eqnarray*}
(1-[-r]+[-1-r]) f(z) + (1-[-s]+[-1-s]) f(y) + (2[-t]+[-2-t]) f(x) & =
& 0\\
(2-[-r]+[-2-r]) f(z) + (2-[-s]+[-2-s]) f(y) + (1-[-t]+[-1-t]) f(x) & =
& 0\\
\end{eqnarray*}
If $r=1$ and $s, t\in\{0, 2\}$, these equations imply that $f(x)=f(y)= 0$,
thus we obtain a contradiction unless all the remaining equations do
not involve $f(z)$. However, this would imply that both points
occurring once in $A$ as well as their sum have $z$-coordinate 1,
which is absurd. The same argument applies if $s=1$ and $r, t\in\{0, 2\}$,
and we deduce the following.

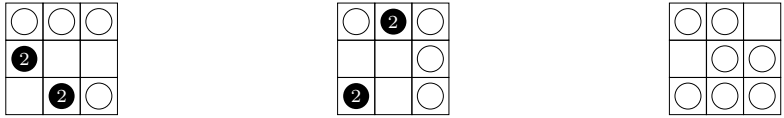
\begin{miniboard}
O O O
2 + +
+ 2 O
\end{miniboard}\hspace*{\fill}
\begin{miniboard}
O 2 O
+ + O
2 + O
\end{miniboard}\hspace*{\fill}
\begin{miniboard}
O O +
+ O O
O O O
\end{miniboard}

Moreover, applying this argument to the sum of the two remaining
elements, we see that this sum cannot have $x$-coordinate $\neq 1$ and
precisely one of $y$ and $z$-coordinate equal to 1.

Suppose that $(0, 1, 1)\in A$. Then this argument gives the following.

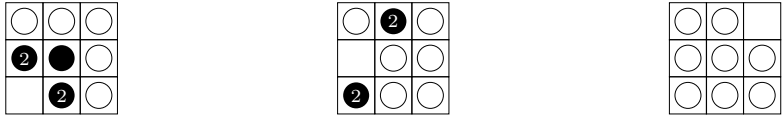
\begin{miniboard}
O O O
2 X O
+ 2 O
\end{miniboard}\hspace*{\fill}
\begin{miniboard}
O 2 O
+ O O
2 O O
\end{miniboard}\hspace*{\fill}
\begin{miniboard}
O O +
O O O
O O O
\end{miniboard}

The sum of all elements equals $(1, 2, 1)$, thus if the remaining
element is in the uppermost layer, there is a zerosum of length
9. Hence, there is only one possibility left, which leads to the matrix
\[
\begin{pmatrix}
1 & 0 & 0\\ 1 & 0 & 0\\ 0 & 1 & 0\\ 0 & 1 & 0\\ 0 & 1 & 0\\ 0 & 1 & 1
\end{pmatrix}
\]
which has rank 3. Hence, $(0, 1, 1)\not\in A$.

Next suppose that $(1, 1, 1)\in A$. Then we have only the following
possibility left.

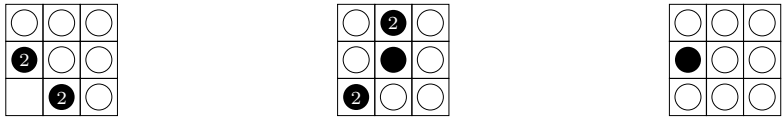
\begin{miniboard}
O O O
2 O O
+ 2 O
\end{miniboard}\hspace*{\fill}
\begin{miniboard}
O 2 O
O X O
2 O O
\end{miniboard}\hspace*{\fill}
\begin{miniboard}
O O O
X O O
O O O
\end{miniboard}

From this we obtain the matrix (deleting the zero-rows coming from
$(1, 1, 1)$) 
\[
\begin{pmatrix}
0 & 1 & 0\\ 1 & 1 & 0\\ 1 & 0 & 1\\ 1 & 0 & 0
\end{pmatrix},
\]
which has rank 3.

Now suppose that $(1, 0, 1)\in A$. Then the only possibility is

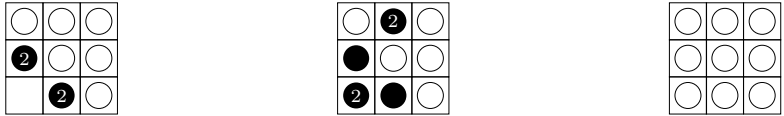
\begin{miniboard}
O O O
2 O O
+ 2 O
\end{miniboard}\hspace*{\fill}
\begin{miniboard}
O 2 O
X O O
2 X O
\end{miniboard}\hspace*{\fill}
\begin{miniboard}
O O O
O O O
O O O
\end{miniboard}

which gives the matrix
\[
\begin{pmatrix}
0 & 1 & 0\\ 0 & 1 & 0\\  0 & 0 & 1\\ 0 & 0 & 1\\ 0 & 0 & 0\\ 1 & 0 & 0
\end{pmatrix}.
\]

It is easy to check that none of the remaining possibilities fits with
$(2, 2, 2)$, thus $(2, 2, 2)\not\in A$.
By now we have reached

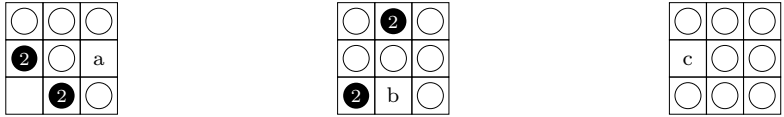
\begin{miniboard}
O O O
2 O a
+ 2 O
\end{miniboard}\hspace*{\fill}
\begin{miniboard}
O 2 O
O O O
2 b O
\end{miniboard}\hspace*{\fill}
\begin{miniboard}
O O O
c O O
O O O
\end{miniboard}

where two of the three cells marekd a, b, c are in $A$. These elements
give the following $2\oplus 3$-matrices:
\[
a\leadsto\begin{pmatrix} 1 & 0 & 0\\ 1 & 1 & 0\end{pmatrix},\quad
b\leadsto\begin{pmatrix} 0 & 0 & 1\\ 0 & 0 & 1\end{pmatrix}, \quad
c\leadsto\begin{pmatrix} 0 & 1 & 0\\ 1 & 1 & 0\end{pmatrix},
\]
obviously, $b$ cannot be in $A$, and the fact that the $x$-coordinate
of $a+c$ is not 1 yields a matrix of rank 3 again. Hence, in each case
we obtain a contradiction.

(viii) We have seen so far that a counterexample to our statement
would have precisely 3 elements occurring twice, which generate
$\Z_3^3$. Hence, we may suppose that $(0, 0, 1)$, $(0, 1, 0)$ and $(1,
0, 0)$ appear twice in $A$, whereas all other elements of $A$ occur
with multiplicity 1. On one hand, having many distinct elements makes
the non-existence of short zerosums a more stringent restriction, on
the other hand, the number of cases increases dramatically. To deal
with this number we use a simple computer program to list all possible
sets of 10 elements without zerosums of length $\leq 3$ or $\geq 8$ having
precisely 3 elements taken twice. With a rather un-sophisticated
approach we obtain a list of 84 configurations, many of which are
symmetric. One easily sees from this list that always one of $(0, 1,
1)$, $(1, 0, 1)$ and $(1, 1, 0)$ is contained in $A$, prescribing this
element by a rotation around the spatial diagonal to be $(0, 1, 1)$
reduces the number of cases to 41, taking care of the remaining
symmetries by hand gives the following list of 16 cases.

\begin{tabular}{cccc}
(0, 1, 1) & (1, 0, 1) & (1, 1, 0) & (1, 1, 1)\\
(0, 1, 1) & (1, 0, 1) & (1, 1, 0) & (1, 1, 2)\\
(0, 1, 1) & (1, 0, 1) & (1, 1, 1) & (1, 1, 2)\\
(0, 1, 1) & (1, 0, 1) & (1, 1, 1) & (1, 2, 0)\\
(0, 1, 1) & (1, 0, 1) & (1, 1, 2) & (1, 2, 2)\\
(0, 1, 1) & (1, 0, 1) & (1, 2, 0) & (1, 2, 1)\\
(0, 1, 1) & (1, 0, 1) & (1, 2, 0) & (1, 2, 2)\\
(0, 1, 1) & (1, 0, 1) & (1, 2, 1) & (1, 2, 2)\\
(0, 1, 1) & (1, 0, 2) & (1, 1, 1) & (1, 1, 2)\\
(0, 1, 1) & (1, 0, 2) & (1, 1, 1) & (1, 2, 1)\\
(0, 1, 1) & (1, 0, 2) & (1, 2, 0) & (1, 2, 1)\\
(0, 1, 1) & (1, 0, 2) & (1, 2, 0) & (1, 2, 2)\\
(0, 1, 1) & (1, 0, 2) & (1, 2, 1) & (1, 2, 2)\\
(0, 1, 1) & (1, 0, 2) & (2, 1, 0) & (2, 1, 1)\\
(0, 1, 1) & (1, 0, 2) & (2, 1, 0) & (2, 1, 2)\\
(0, 1, 1) & (1, 0, 2) & (2, 1, 1) & (2, 1, 2)\\
\end{tabular}

Let $(r, s, t)$ be an element of $A$ occurring once, distinct from
$(0, 1, 1)$. Then we have the zerosums $(0, 1, 1)+2x+2y$, $(r, s,
t)+[-t]x+[-s]y+[-t]z$ and $(0, 1, 1) + (r, s,
t)+[-1-t]x+[-1-s]y+[-1-t]z$, which together imply the equations
\[
\begin{array}{rcrcrcrcrcl}
f((0, 1, 1)) && && & + & 2 f(y) & + & 2 f(x) & = & 1\\
 && f((r, s, t)) & + & [-r] f(z) & + & [-s] f(y) & + & [-t] f(x) & = & 1\\
f((0, 1, 1)) & + & f((r, s, t)) & + & [-1-r] f(z) & + & [-1-s] f(y) & + & [-1-t] f(x) & = & 1
\end{array}
\]
Substracting the third equation from the sum of the other two
equations we obtain
\[
(2+[-s]-[-1-s]) f(y) + (2+[-t]-[-1-t]) f(x) = 1.
\]
Note that for $a\in\Z_3^3$ we have
\[
\frac{2+[-a]-[-1-a]}{3} = \begin{cases} 0, & a=0\\ 1, & a=1, 2\end{cases}.
\]
Now consider the first entry in the table above. The second element is
$(1, 0, 1)$, which gives the equation $f(x) = 1$, the third element
yields $f(y)=1$, whereas the third one implies $f(x)+f(y)=1$, and we
obtain a contradiction. In the same way we can deal with all cases in
which the third an fourth entry contains an entry 0, which is the case
for all but the following.

\begin{tabular}{cccc}
(0, 1, 1) & (1, 0, 1) & (1, 1, 1) & (1, 1, 2)\\
(0, 1, 1) & (1, 0, 1) & (1, 1, 2) & (1, 2, 2)\\
(0, 1, 1) & (1, 0, 1) & (1, 2, 1) & (1, 2, 2)\\
(0, 1, 1) & (1, 0, 2) & (1, 1, 1) & (1, 1, 2)\\
(0, 1, 1) & (1, 0, 2) & (1, 1, 1) & (1, 2, 1)\\
(0, 1, 1) & (1, 0, 2) & (1, 2, 1) & (1, 2, 2)\\
(0, 1, 1) & (1, 0, 2) & (2, 1, 1) & (2, 1, 2)\\
\end{tabular}

Finally, there is no need to take $(r, s, t)$ to be one of the three
elements occurring once distinct from $(0, 1, 1)$, we could also take
the sum of two distinct ones among them. Since in the 7 remaining
cases we already have found equations coming from $s=0, t\neq 0$ and $s,
t\neq 0$, it suffices to obtain $(r, s, t)$ with $s\neq 0, t=0$. For the
first three cases this is achieved by adding the second element to the
fourth one, whereas for the other 4 cases we add the second element to
the third one. Hence, in all cases we reach a contradiction, which
finishes our proof.
\end{proof}

From Theorem~\ref{thm:Af} we can now deduce Theorem~\ref{thm:D=MZ3}.

\begin{proof}[Proof of Theorem~\ref{thm:D=MZ3}.]
We may suppose that $(6, d)=1$, since for all other cases this is
already known.
Let $A$ be a set of $3d+4\geq 13$ points such that there is a function as
above. Then $A$ contains a zerosum of length $\leq 3$. In fact, the
argument at the beginning of the proof of the previous theorem did not
involve long zero-sums and yields that in any set without zerosums of
length 3 admitting a function as described there cannot be 5 elements
occurring twice. Hence, there are at least 9 distinct elements in $A$,
but we know that among 9 distinct elements there is always a zerosum
of length $\leq 3$. Thus, if there is a set $B$ in $G$ of size $3d+4$
without a zerosum, there exists a set $A$ in $\Z_3^3$ of the same size
with at most $d-1$ disjoint zerosum subsets admitting a function $f$
as above. Clearly, the function $f$ is also admissible for every subset
of $A$, and we can remove $d-2$ zerosum subsets of size $\leq 3$ to reach
a set of at least 10 points admitting a function, which does not have
2 disjoint zerosum subsets. However, by the previous theorem such a function
does not exist. Hence, there does not exist a zerosum free set of size
$3d+4$, and we deduce $D(G)\leq 3d+4$. The other inequality is trivial.
\end{proof}

\begin{tabular}{ll}
Gautami Bhowmik, & Jan-Christoph Schlage-Puchta,\\
Universit\'e de Lille 1, & Albert-Ludwigs-Universit\"at,\\
Laboratoire Paul Painlev\'e, & Mathematisches Institut,\\
U.M.R. CNRS 8524,  & Eckerstr. 1,\\
  59655 Villeneuve d'Ascq Cedex, & 79104 Freiburg,\\
  France & Germany\\
bhowmik@math.univ-lille1.fr & jcp@math.uni-freiburg.de
\end{tabular}
\end{document}

%% file: Dav.board.tex
\usepackage{amsmath,amssymb,multicol}
\usepackage[active]{srcltx}

\newcommand{\Z}{\mathbb{Z}}

\newcommand{\fieldsize}{14} 
\newcommand{\fieldsizebp}{\fieldsize bp}

\newcommand{\polycellsize}{6} 

\newcommand{\fieldend}{%
\hskip\fieldsizebp\vrule width0mm height\fieldsizebp\relax%
}

\newcommand\specialps{} 
\begingroup
\catcode`\"=12 
\gdef\specialps#1{\special{" field-setup #1}}
\endgroup

\newcommand{\fieldps}[1]{%
  \vrule width0mm height\fieldsizebp\relax
  \specialps{#1}%
}

\newcommand{\fieldcenter}[1]{\rlap{\vbox to\fieldsizebp{\vfil\hbox to\fieldsizebp{{\rm\hfil \strut #1\hfil}}\vfil}}}
\newcommand{\fieldcenterwhite}[1]{
\fieldcenter{#1}
}

\newcommand{\extN}{}
\newcommand{\extS}{}
\newcommand{\extE}{}
\newcommand{\extW}{}
\newcommand{\fif}{ff \extN\extS\extE\extW}

\newcommand{\nostone}{\fieldps{\fif}\fieldend}

\newcommand{\nostoned}{\fieldps{\fif do}\fieldend}

\newcommand{\fieldstonew}[1]{\fieldps{\fif sw}\fieldcenter{\scriptsize #1}\fieldend}
\newcommand{\fieldstoneb}[1]{\fieldps{\fif sb}\fieldcenterwhite{\scriptsize #1}\fieldend}

\newcommand{\nothing}{\fieldend}

\newcommand{\letter}[1]{\fieldps{\fif}\fieldcenter{{\footnotesize #1}}\fieldend}

\newcommand{\fieldministonew}{\fieldps{\fif mw}\fieldend}
\newcommand{\fieldministoneb}{\fieldps{\fif mb}\fieldend}
\newcommand{\fieldministoned}{\fieldps{\fif md}\fieldend}

\newcommand{\paveup}{\fieldps{\fif pu}\fieldend}
\newcommand{\pavedown}{\fieldps{\fif pd}\fieldend}
\newcommand{\paveleft}{\fieldps{\fif pl}\fieldend}
\newcommand{\paveright}{\fieldps{\fif pr}\fieldend}
\newcommand{\paveupdown}{\fieldps{\fif pud}\fieldend}
\newcommand{\paveleftright}{\fieldps{\fif plr}\fieldend}
\newcommand{\flagstone}{\fieldps{\fif fs}\fieldend}

\newcommand{\fieldcoord}[1]{\vrule width0mm height\fieldsizebp\relax\fieldcenter{#1}\fieldend}

\makeatletter

\newcommand{\boardcatcodes}{%
\catcode`\^^M\active
\catcode`\ \active
\catcode`\.\active
\catcode`\O\active
\catcode`\X\active
\catcode`\0\active
\catcode`\1\active
\catcode`\2\active
\catcode`\3\active
\catcode`\4\active
\catcode`\5\active
\catcode`\6\active
\catcode`\7\active
\catcode`\8\active
\catcode`\9\active
\catcode`\*\active
\catcode`\+\active
\catcode`\o\active
\catcode`\x\active
\catcode`\q\active
\catcode`\a\active
\catcode`\b\active
\catcode`\c\active
\catcode`\d\active
\catcode`\e\active
\catcode`\f\active
\catcode`\g\active
\catcode`\h\active
\catcode`\>\active
\catcode`\<\active
\catcode`\^\active
\catcode`\v\active
\catcode`\|\active
\catcode`\-\active
\catcode`\#\active
\catcode`\\\active
\catcode`\N\active
\catcode`\S\active
\catcode`\E\active
\catcode`\W\active
}

\def\definitionstack{}

\long\def\addboarddefinition#1#2{%
\expandafter\def\expandafter\definitionstack\expandafter{\definitionstack%
\def#1{#2}%
}%
}

\def\boarddefinition{%
\begingroup%
\boardcatcodes%
\boarddefinition@%
}
\def\boarddefinition@#1{%
\endgroup%
\addboarddefinition{#1}
}

\boarddefinition{N}{\def\extN{N }}
\boarddefinition{S}{\def\extS{S }}
\boarddefinition{E}{\def\extE{E }}
\boarddefinition{W}{\def\extW{W }}
\boarddefinition{ }{\doline\doboardspace}
\boarddefinition{+}{\doline\nostone\eatspacemode}
\boarddefinition{.}{\doline\nostoned\eatspacemode}
\boarddefinition{*}{\doline\nostoned\eatspacemode}
\boarddefinition{O}{\doline\stonemodes{\fieldstonew}}
\boarddefinition{X}{\doline\stonemodes{\fieldstoneb}}
\boarddefinition{o}{\doline\fieldministonew\eatspacemode}
\boarddefinition{x}{\doline\fieldministoneb\eatspacemode}
\boarddefinition{q}{\doline\fieldministoned\eatspacemode}
\boarddefinition{
}{\donewline} 
\boarddefinition{0}{\doline\savenumbereven\doboardnumber{0}}
\boarddefinition{1}{\doline\savenumberodd\doboardnumber{1}}
\boarddefinition{2}{\doline\savenumbereven\doboardnumber{2}}
\boarddefinition{3}{\doline\savenumberodd\doboardnumber{3}}
\boarddefinition{4}{\doline\savenumbereven\doboardnumber{4}}
\boarddefinition{5}{\doline\savenumberodd\doboardnumber{5}}
\boarddefinition{6}{\doline\savenumbereven\doboardnumber{6}}
\boarddefinition{7}{\doline\savenumberodd\doboardnumber{7}}
\boarddefinition{8}{\doline\savenumbereven\doboardnumber{8}}
\boarddefinition{9}{\doline\savenumberodd\doboardnumber{9}}
\boarddefinition{a}{\doline\doboardletter{a}}
\boarddefinition{b}{\doline\doboardletter{b}}
\boarddefinition{c}{\doline\doboardletter{c}}
\boarddefinition{d}{\doline\doboardletter{d}}
\boarddefinition{e}{\doline\doboardletter{e}}
\boarddefinition{f}{\doline\doboardletter{f}}
\boarddefinition{g}{\doline\doboardletter{g}}
\boarddefinition{h}{\doline\doboardletter{h}}
\boarddefinition{>}{\doline\paveright\eatspacemode}
\boarddefinition{<}{\doline\paveleft\eatspacemode}
\boarddefinition{^}{\doline\paveup\eatspacemode}
\boarddefinition{v}{\doline\pavedown\eatspacemode}
\boarddefinition{|}{\doline\paveupdown\eatspacemode}
\boarddefinition{-}{\doline\paveleftright\eatspacemode}
\boarddefinition{#}{\doline\flagstone\eatspacemode}
\boarddefinition{\}{\dobackslash}
{{{{{{}}}}}}  

\newcommand{\activateboard}{%
\begingroup
\definitionstack
\boardcatcodes
}

\newcommand{\deactivateboard}{%
\coordx\endgroup%
}

\def\activateboarddefs{%
\def\normalnewlinemode{\def\donewline{\doline\doboardspace\par\normalspacemode\beginlinemode}}
\def\eatnewlinemode{\def\donewline{\normalnewlinemode}}
\eatnewlinemode
\def\normalspacemode{\def\doboardspace{\nothing\eatspacemode}}%
\def\specialspacemode##1{%
  \def\doboardspace{##1\normalmodes}%
}%
\def\eatspacemode{\specialspacemode{}}%
\def\savenumberodd{\def\reallydonumber####1{\fieldstoneb{####1}}}%
\def\savenumbereven{\def\reallydonumber####1{\fieldstoneb{####1}}}%
\def\normalnumbermode{\def\doboardnumber####1{%
   \def\savenumber{####1}%
   \specialspacemode{\reallydonumber{\savenumber}}%
   \specialnumbermode{\reallydonumber{\savenumber\extranumber}}%
}}%
\def\specialnumbermode##1{%
  \def\doboardnumber####1{\def\extranumber{####1}##1\normalmodes}%
}
\def\normallettermode{\def\doboardletter####1{%
  \def\saveletter{####1}%
  \specialspacemode{\letter{\saveletter}}%
  \specialnumbermode{\letter{\saveletter$_{\extranumber}$}}%
}}%
\def\speciallettermode##1{%
  \def\doboardletter####1{\def\extraletter{####1}##1\normalmodes}%
}
\def\stonemodes##1{%
  \specialspacemode{##1{}}%
  \speciallettermode{##1{\extraletter}}%
}%
\def\normalmodes{\normalspacemode\normalnumbermode\normallettermode}
\normalmodes
\def\dobackslash{\deactivateboard\@ifnextchar e{\readend}{\readnotes}}
\def\readend end{\end}
\def\readnotes notes{%
\endgroup 
\begingroup 
\parindent0mm
\small}
\def\beginlinemode{\def\doline{\coordy\normallinemode}}%
\def\normallinemode{\def\doline{}}%
\beginlinemode%
\def\nocoordinates{%
  \def\coordx{}%
  \def\coordy{}%
  \activateboard%
}%
\def\readcoordinates(##1,##2){%
  \def\coordx{%
    \nothing%
    \setcounter{coordcnt}{1}%
    \loop\fieldcoord{\alph{coordcnt}}\ifnum\c@coordcnt<##1\addtocounter{coordcnt}{1}\repeat%
    \par
  }%
  \def\coordy{%
    \fieldcoord{$\arabic{coordcnt}$}%
    \addtocounter{coordcnt}{-1}%
  }%
  \setcounter{coordcnt}{##2}%
  \activateboard%
}%
\def\readmorecoordinates(##1,##2)(##3,##4){%
  \def\coordx{%
    \nothing%
    \setcounter{coordcnt}{##1}%
    \loop\fieldcoord{{\small$\arabic{coordcnt}$}}\ifnum\c@coordcnt<##3\addtocounter{coordcnt}{1}\repeat%
    \par
  }%
  \def\coordy{%
    \fieldcoord{{\small$\arabic{coordcnt}$}}%
    \addtocounter{coordcnt}{-1}%
  }%
  \setcounter{coordcnt}{##4}%
  \activateboard%
}%
}

\makeatother

\newcommand{\boardwidth}{0.43\textwidth}

\newenvironment{bboard}{%
\begin{minipage}{\boardwidth}%
\vskip1ex%
\begingroup%
\rm%
\parindent0mm%
\baselineskip0mm%
\lineskiplimit0mm%
\lineskip0mm%
\activateboarddefs%
}{%
\endgroup%
\vspace*{1ex}%
\end{minipage}\hspace{\fill}%
}

\newenvironment{miniboard}{%
  \begingroup\renewcommand{\boardwidth}{0.2\textwidth}%
  \begin{bboard}\nocoordinates%
}{\end{bboard}\endgroup}

\newcounter{coordcnt}

\newlength{\polycelllength}
\setlength{\polycelllength}{\polycellsize bp}

\newcommand{\polyo}[2]{\strut\specialps{#1}\kern#2\polycelllength}

%% file: Dav3.bbl
\begin{thebibliography}{99}
\bibitem{Boas1} P. van Emde Boas, D. Kruyswijk, A combinatorial problem on
  finite abelian groups III, Rapport ZW-1969-008, Mathematisch Centrum Amsterdam.
\bibitem{DOQ} C. Delorme, O. Ordaz, D. Quiroz, Some remarks on
  Davenport constant, {\em Discrete Math.} {\bf 237} (2001), 119--128.
\bibitem{Olson} J. E. Olson, An addition theorem modulo $p$, {\em 
J. Combinatorial Theory} {\bf 5} (1968) 45--52.
\end{thebibliography}
